\documentclass[preprint,12pt]{elsarticle}

\usepackage{amssymb}
\usepackage{amsmath}

\newcounter{bla}

\journal{journal for publication}

\usepackage{graphicx}
\usepackage[normalem]{ulem}
\usepackage{color}
\usepackage{dcolumn}
\usepackage{bm}
\usepackage{float}

\begin{document}

\begin{frontmatter}



\title{Time-Spectral Efficiency}


\author{Jan Scheffel \corref{author}}

\cortext[author] {Jan Scheffel.\\\textit{E-mail address:} jans@kth.se}
\address{Electromagnetic Engineering and Fusion Science \\ 
	KTH Royal Institute of Technology, Stockholm, Sweden}

\begin{abstract}
This study concerns the efficiency of time-spectral methods for numerical solution of differential equations. It is found that the time-spectral method GWRM demonstrates insensitivity to stiffness and chaoticity due to the implicit nature of the solution algorithm. Accuracy is thus determined primarily by numerical resolution of the solution shape. Examples of efficient solution of stiff and chaotic problems, where explicit methods fail or are significantly slower, are given. Non-smooth and partially steep solutions, however, remain challenging for convergence and accuracy. Some, earlier suggested, smoothing algorithms are shown to be ineffective in addressing this issue. Our findings underscore the need for further exploration of time-spectral approaches to enhance convergence and accuracy for steep or non-smooth solutions.
\end{abstract}

\begin{keyword}
Time-spectral; GWRM; stiffness; chaos; stiff; chaotic; ODE; PDE; RK4; numerical efficiency.

\end{keyword}

\end{frontmatter}

\section{Introduction}
\noindent Time-spectral methods have been proposed as efficient alternatives to explicit and implicit time-stepping methods for solving initial-value differential equations \cite{Scheffel:GWRM1,Scheffel:GWRM2}. A number of examples and applications, including solution of chaotic problems related to numerical weather prediction \cite{Scheffel:GWRM3}, 1D problems in resistive magnetohydrodynamics \cite{Scheffel:GWRM5} and 2D nonlinear resistive magnetohydrodynamic turbulence \cite{Scheffel:GWRM6} can be found in the literature. In this article, we will generally and critically assess time-spectral methods with regards to accuracy and efficiency. 

When evaluating the efficiency of a time-spectral method, or any computational approach for differential equations, at least four properties of the problems addressed come into play: 1) \textit{steepness} (normalized maximum rate of change), 2) \textit{oscillations} (number of extrema in the computational interval), 3) \textit{stiffness} (at least one Lyapunov exponent $\gamma_i < 0$, with large magnitude) and 4) \textit{chaoticity} (at least one $\gamma_i > 0$). We will consider these aspects, making use of the concept of local Lyapunov exponents $\gamma_i $ (LLEs), as described in \cite{Scheffel_transforming, Cartwright}, for diagnosing stability of the solution in local sub-intervals of the temporal domain. Asymptotically stable solutions are characterized by $\gamma_i < 0$, with $|\gamma_i (\Delta T)^{-1}| \lesssim 1$, where $\Delta T$ is the sub-interval or time step length, and correspond to non-stiff and non-chaotic problems.

The time-spectral Generalized Weighted Residual Method (GWRM) is briefly outlined in Section 2. This initial-value method is based on a solution ansatz in terms of truncated expansions of Chebyshev polynomials in time and space, and also in parameter space when so desired. In the latter case the parameter dependence of a solution can be found from a single computation. The coefficients of the ansatz are obtained from a global weighted average method, resulting in a system of algebraic equations, to be solved iteratively for nonlinear problems. Since the final solution is an analytical function of time and, for PDEs, additional space variables, it conveniently lends itself to further analytical manipulation. The time intervals employed are usually one to two orders of magnitude longer than the time steps of finite difference methods. The GWRM has been further developed for a number of applications. This includes designing the iterative, algebraic equation solver SIR, which features high global convergence \cite{Scheffel:SIR,Scheffel:SIR2}, as well as subdomain methods \cite{Scheffel:GWRM11,Scheffel:GWRM21} and automatic time interval adaption \cite{Scheffel:GWRM10}.

In Section 3, two recently proposed methods for smoothening in problems with steep or oscillating solutions will be briefly summarized. These are shown to be well adapted to time-spectral methods.

The concept of Lyapunov exponents is helpful when identifying stiffness and chaoticity in nonlinear problems. A novel scheme for approximative computation of \textit{localized} Lyapunov exponents (LLEs), without needing to solve the differential equations \cite{Scheffel_transforming}, is employed in Section 4. Two problems, one stiff and one chaotic, are specified as aids for evaluating GWRM efficiency in Section 5. Here, an overall assessment of GWRM efficiency with regards to convergence, accuracy and efficiency in relation to the problem properties of steepness, oscillations, stiffness and chaoticity will be made. Comparisons to explicit and implicit time-stepping methods will also be performed.

\section{The Generalized Weighted Residual Method (GWRM)}

\noindent The GWRM \cite{Scheffel:GWRM1,Scheffel:GWRM2} is a \textit{time-spectral method}, thus differing from traditional spectral methods for initial-value problems \cite{Gottlieb:1} in that also the time domain is given a spectral representation. Consider the following system of initial-value PDEs
\begin{equation}
\frac{\partial\mathbf{u}}{\partial t}=D\mathbf{u}+\mathbf{f}
\end{equation}

Here $D$ denotes a linear or nonlinear matrix differential operator that may depend on physical variables ($t$, $\mathbf{x}$, and $\mathbf{u}$) as well as on physical parameters (denoted $\mathbf{p}$), and $\mathbf{f}={\mathbf{f}}(t,\mathbf{x};\mathbf{p})$ is a known source or forcing term. An approximate solution ansatz is given as a truncated multivariate series of first kind Chebyshev polynomials $T_n(z)$;
\begin{equation}
u(t,x;p)=\sum_{k=0}^{K}{'} \sum_{l=0}^{L}{'}\sum_{m=0}^{M}{'}a_{klm}T_k(\tau)T_l(\xi)T_m(P)
\end{equation}
for one spatial dimension $x$ and one parameter $p$. We have defined

\begin{equation}
\tau \equiv \frac{t-A_t}{B_t},
\xi \equiv \frac{x-A_x}{B_x},
P \equiv \frac{p-A_p}{B_p},
\end{equation}
\begin{equation}
A_z \equiv \frac{z_1+z_0}{2},
B_z \equiv \frac{z_1-z_0}{2}
\end{equation}

\noindent where $z$ can be any of $t$, $x$ or $p$ and $z_0$, $z_1$ are values at interval boundaries. Primes denote that 0th order coefficients should be halved. The coefficients $a_{klm}$ of this \textit{analytical} expression obtains from weighted residuals as follows. For a single differential equation in $u$, a residual $R$ is defined as \cite{Fletcher:1}: 
\begin{equation}
R=u(t,{x};{p})-\left(u(t_0,{x};{p})+\int_{t_0}^{t}(Du+f)dt'\right)
\end{equation}

An integration is performed over the entire computational domain:

\begin{equation}
\int_{t_0}^{t_1}\int_{{x_0}}^{{x_1}}\int_{{p_0}}^{{p_1}}RT_q(\tau)T_r(\xi)T_s(P)w_tw_xw_pdtd{x}d{p}=0  .
\end{equation}

The following weight factors are used \cite{Mason:1}:
\begin{equation}
w_t=(1-\tau^2)^{-\frac{1}{2}},
w_x=(1-\xi^2)^{-\frac{1}{2}},
w_p=(1-P^2)^{-\frac{1}{2}}
\end{equation}

The Chebyshev coefficients $b_{lm}$ of the initial condition, given by $u(0,x;p)$, should now be obtained, as well as the coefficients $A_{klm}$ from the expansion
\begin{equation}
\int_{t_0}^{t}Dudt'=\sum_{k=0}^{K}{'}\sum_{l=0}^{L}{'}\sum_{m=0}^{M}{'}A_{klm}T_k(\tau)T_l(\xi)T_m(P)  ,
\end{equation}
showing that Eq. (8) is a linear/nonlinear equation in the parameters $a_{klm}$ depending on the form of Eq. (1). The coefficients $F_{klm}$ (relating to $\bf f$), are obtained by a similar procedure to that in Eq. (8). 

There results the following simple system of algebraic equations for the Chebyshev coefficients of the solution (2): 
\begin{equation}
	a_{klm}=2\delta_{k0}b_{lm}+A_{klm}+F_{klm}
\end{equation}
Boundary conditions enter by replacing high-end $l$-indexed coefficients of Eq. (9) with coefficients relating to corresponding boundary equations. 

A semi-implicit root solver (SIR) has been developed \cite{Scheffel:SIR,Scheffel:SIR2} to iteratively and efficiently solve the system of implicit non-linear algebraic equations for the Chebyshev coefficients in Eq. (9). SIR generally has superior convergence properties to Newton's method. Furthermore the algorithm avoids landing on local minima, as may occur for Newton's method with line-search. 

Two features of the Standard GWRM (SGWRM), as described above, are worthy of consideration in the present context. First, there exist optimal, local interval lengths. It was shown in \cite{Scheffel:GWRM3} that whereas employing large values of $K,L$ or $M$ enables high accuracy in long intervals, the overhead in numerical computations (determining $A_{klm}$, for example) becomes costly. Conversely, for small values of $K,L$ or $M$ and many shorter intervals, even in combination with less iterations, the number of calls for the algorithm that computes Eq. (8) becomes prohibitive. The reason Eq. (8) usually needs continual symbolic recomputation is that adaptive time interval control generally is favourable.  

Second, although the GWRM is designed to solve causal initial-value problems, \textit{the numerical algorithm is acausal} in the sense that Eq. (9) is simultaneously iterated implicitly for all coefficients $a_{klm}$ of the solution $u(t)$ for the entire temporal interval. This is similar to the handling of spatial boundaries in PDE problems. As a result, CFL-like causal restrictions on time interval length are eliminated, as well as numerical sensitivity to local, rapidly changing spurious solutions. This property is strongly favourable when addressing stiff and chaotic problems, as shown in Sections 4 and 5.

\section{Steepness and oscillations - temporal smoothing}

\noindent Numerical tools for smoothening in differential equation problems with oscillative and locally steep solutions have been introduced elsewhere, with applications to example problems \cite{Scheffel:GWRM7}. Applications to the GWRM, also for stiff and chaotic problems, will be critically studied here. The discussion will be confined to ODEs only, but the results can straightforwardly be generalized to PDEs. The following measure of steepness, or rather non-smooth functional behaviour of $u(t)$ in an interval $[0,T]$, was put forth:
\begin{equation}
	S_u \equiv \underset{t \in [0,T]}{max} \frac{\lvert du/dt \rvert}{(u_{max}-u_{min})/T}
\end{equation}

\noindent Here $u_{max}$ and $u_{min}$ are the maximum and minimum values of $u(t)$ in $[0,T]$. With this definition, the smoothest function is a straight line, hence $S_u\geq 1$. Since time-spectral methods are inefficient for strongly non-smooth problems, requiring a large number of Chebyshev modes for accurate resolution, methods for solution smoothing are here of great interest.

\subsection{The Time Averaging (TA) Method}

\noindent A natural choice for a smoothing operator is a symmetric \textit{running average} of the solution $u(t)$:
\begin{equation}
	U(t)\equiv \frac{1}{2\Delta}\int_{t-\Delta}^{t+\Delta} u(t') dt'
\end{equation}

\noindent where the interval length $2\Delta$ is arbitrary. Sometimes the equivalent definition

\begin{equation}
	U(t)\equiv \frac{1}{2\Delta}\int_{-\Delta}^{\Delta} u(t+t') dt'
\end{equation}

\noindent is employed \cite{Acharya2006:1}. It was found in \cite{Scheffel:GWRM7} that general procedures for transforming ODE's to \textit{exact} differential equations for $U(t)$ obtains by using the relation

\begin{equation}
	\frac{dU}{dt} =  \frac{1}{2\Delta} (u(t+\Delta)-u(t-\Delta))
\end{equation}

For linear ODEs with constant coefficients \textit{the only change} in the problem formulation is a modification of the initial condition. The running average of the solution is thus determined for the entire solution interval by a simple alteration of the initial condition. 

For nonlinear ODEs, or linear ODEs with variable coefficients, the number of ODEs \textit{doubles} since a linearly independent function of $u(t+\Delta)$ and $u(t-\Delta)$ must also be introduced. The \textit{running two-point average} is chosen:
\begin{equation}
	V^{(n)}(t) \equiv \frac{1}{2}(u(t+\Delta)+u(t-\Delta))
\end{equation}

\noindent where $V^{(n)}(t) \equiv d^nV(t)/dt^n$ and $V^{(0)}(t)$ stands for $V(t)$. Employing Eqs. (13) and (14) all ODEs could be rewritten as exact ODEs for the two running averages. By choosing $n \geq 1$ and solving for $V(t)$ both $U(t)$ and $V(t)$ may be expected to be relatively smooth, since they are time integrals of combinations of the sought solution. Thus less temporal basis functions are needed for the same accuracy as for solution of the original problem. This approach is termed the Time Averaging (TA) Method.

\subsection{Long Time Averaging (LTA) and Time Integration (TI) Methods}

\noindent It is of interest to, alternatively, solve for an averaged solution that retains the original number of ODEs. This is accomplished by redefining the ODEs as a set of exact equations for the \textit{long-time average}
\begin{equation}
	W(t)\equiv \frac{1}{t}\int_{0}^{t} u(t') dt' .
\end{equation}

The recipe is to substitute $u(t)=d(tW)/dt \equiv dZ/dt$ into the original set of ODEs, solve for $W(t)$ and finally either settle for $W(t)$ or transform back to $u(t)$. For the model equation

\begin{equation}
	\frac{du(t)}{dt} = F(t,u(t)),  \quad u(0) = u_0
\end{equation}

\noindent where $F$ is a linear or nonlinear algebraic operator, there results the following integro-differential equation to be solved:
\begin{equation}
	Z(t) = u_0t+\int_{0}^{t} \int_{0}^{t'} F(t'',dZ/dt'') dt''dt'
\end{equation}

\noindent where $Z(0)=0$ without loss of generality. Finally, $W(t)=Z(t)/t$ and $u(t)=dZ/dt$ obtains. This formulation, the Long-Time Averaging (LTA) Method, is well suited for the GWRM. 

The Time Integration (TI) Method \cite{Scheffel:GWRM7}, is strongly related to the LTA Method. In order to simultaneously a) reformulate the set of ODEs to a problem for a smoother solution function and b) retain the number of ODEs, the function
\begin{equation}
	v(t)\equiv \int_{0}^{t} u(t') dt' + At .
\end{equation} 

\noindent is introduced into the ODE, using the relation $dv/dt=u(t)+A$. It is expected that the solution $v(t)$ is smoother than $u(t)$ itself. The constant $A$ is motivated by that a suitable choice will render the steepness \cite{Scheffel:GWRM7} $S_v$ (Eq. (10)) of $v(t)$ in the interval $[0,T]$ minimal in some sense.

It is readily seen that the TI method is equivalent to the LTA method if a constant $A$ is added to the right hand side of Eq. (15). The relation $W(t)=v(t)/t$ is then obtained. Since $v(t)$ is the computed analytical function in the TI Method, the long-term average $W(t)$ is immediately obtained analytically, and there is little reason to compute $W(t)$ at all. Thus the LTA Method is fully contained in the TI Method.

\subsection{Resolving oscillations}
\noindent It is of interest to estimate the number of modes $K_a=K_a(N_{e},\epsilon)$ required for approximating a function $f(t)$, featuring $N_{e}$ extrema in an interval $[t,t+\tau]$, within a relative error $\epsilon$. This question is discussed in \cite{Gottlieb:1} where it is argued that, for \textit{smooth oscillating functions}, Chebyshev series expansions with approximately $\pi$ polynomials per "wavelength" should be employed for an accuracy of about one percent. By studying a set of oscillating functions, we can qualify this statement somewhat. It is found that, for reasonably periodical functions, maximum relative errors $\epsilon$ of 0.01, 0.001, and 0.0001 were obtained by employing approximately 2, 2.5 and 3 Chebyshev modes respectively for each extremum of $f(t)$. 

For Chebyshev approximation of \textit{aperiodically oscillating functions} $f(t)$, the following approximate relations, valid for intervals in which $N_e \in [1,6]$, were found from a statistical set using randomly generated amplitudes and phases: 
	\begin{align}
	K_a(N_{e},0.01) = 1.5 N_e + 3.5  \\
	K_a(N_{e},0.001) = 1.7 N_e + 4.4
	\end{align}

GWRM solution of ODEs and PDEs, however, generally require \textit{additional modes} for the same accuracy, due to the need for accurate representation of the differentiated functions. The dominant error typically occurs near the temporal interval boundaries. For \textit{temporal system order} $O_t$, the term $O_t$ should be added to the right hand side of Eqs. (19) and (20).

\section{Stiff and chaotic nonlinear problems}
\noindent We will here present two nonlinear problems, one stiff and one chaotic, that will be employed in the next Section for general assessment of the performance of the GWRM in comparison with time-stepping methods. \textit{Stiffness} appears when neighbouring solutions feature a very rapid time dependence, so that the temporal resolution for solving the problem, using explicit methods, is much higher than what is required by the curvature of the solution trajectory. \textit{Chaoticity} occurs for differential equations, where neighbouring solutions to the exact solution strongly diverges in time, a phenomenon which is characterized by positive (local or non-local) Lyapunov exponents.   

\subsection{Autocatalytical chemical reactions - a stiff problem}
\noindent The following system of equations emerge in a kinetic analysis of chemical autocatalytic reactions \cite{Robertson}:
\begin{equation}
	\begin{split}
	\frac{dx}{dt} = -ax+byz  \\
	\frac{dy}{dt} = ax-byz-cy^2  \\
	\frac{dz}{dt} = cy^2 
	\end{split}
\end{equation}

\noindent where $a=0.04$, $b=10^4$, and $c=3\cdot10^7$. Initial conditions are $(x,y,z)=(1,0,0)$. Since the time scales of the reactions are very different, due to the large difference between the (reaction rate) constants, explicit methods are likely to fail for long intervals like $[10^{-6},10^6]$. This is an example of stiffness. As can be seen from the time-adaptive SGWRM run in Fig. 1, $y(t)$ features a strong transient initially (logarithmic time axis), whereafter all three variables $x(t),y(t)$ and $z(t)$ feature smooth solutions. The accuracy of the root solver SIR was set to $10^{-6}$. The CPU time was 1.2 s and 40 Mb of memory was used. 

\begin{figure}[h!]
	\centering
	\includegraphics[width=5in]{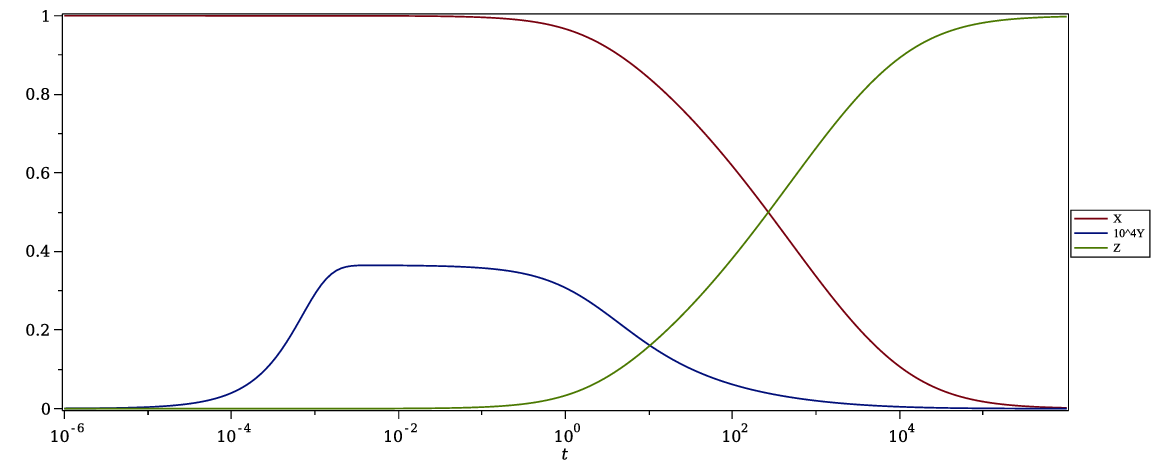}
	\caption{SGWRM time-adaptive solution of Robertson equations, using 49 time intervals for $K=6$, with accuracy $\epsilon=0.001$. Note that $y(t)$ has been magnified a factor $10^4$.}
\end{figure}

Considering that this problem is very stiff, we may anticipate difficulties when attempting to solve this problem using explicit time-stepping methods. 

The stiffness of Eq. (21) was studied using local Lyapunov exponents (LLEs) $\gamma_i$ in \cite{Scheffel_transforming}. For an interval $[t,t+\tau]$, these are defined as

\begin{equation}
	\gamma_i(\tau,t) =  \,  \lim_{|\delta X_{i}| \rightarrow 0} \, \frac{1}{\tau} ln \frac{|\delta X_i(t+\tau)|}{|\delta X_i(t)|}  \quad ,
\end{equation}   

\noindent where $\delta X_i$ represents the deviations of two neighbouring solutions in the variable $X_i$ in $N$-dimensional phase space. Stiff differential equations feature negative $\gamma_i$ with large magnitudes, whereas chaoticity holds for positive $\gamma_i$.    

At the initial state, the Robertson equations are not stiff; here $\gamma_1=0$, $\gamma_2=0$ and $\gamma_3=-0.04$. But already at $y=10^{-6}$ (for $x=1$ and $z=0$) the system becomes stiff; $\gamma_1=0$, $\gamma_2=-0.05$ and $\gamma_3=-60$. Extreme stiffness is obtained for $t>10^{-4}$, as $y>0.4\cdot10^{-4}$ and $\gamma_3 \le -2400$. 

A fourth order \textit{explicit} Runge-Kutta (RK4) algorithm with adaptive step size was employed to solve the problem (21). After 500 CPU-seconds, however, the algorithm stagnated near $t=96$. Over 61 000 time steps were used at this stage in the computation. Similar difficulties were encountered for other parameters. It is clear that explicit time-stepping methods are of limited use for stiff problems of this type.

Expecting \textit{implicit} time-stepping algorithms to be a better choice for stiff problems, a standard trapezoid algorithm with adaptive step size was also applied. For an accuracy $\epsilon=0.001$, a solution was obtained after 13 CPU seconds, using 118 steps for the whole interval with initial step length 0.1. A similar step size adaption algorithm as for the GWRM solution was used. The implicit trapezoid algorithm was thus nearly 10 times less efficient than the GWRM.

\subsection{The Lorenz 1984 equations - a chaotic problem}
\noindent In 1984, Lorenz \cite{Lorenz:1984atmosphere} formulated a set of three non-linear ordinary differential equations that features chaotic behaviour similar to that characterizing meteorological systems. It is a simple model of what is termed Hadley circulation. Albeit the Lorenz 1984 model is not an accurate NWP model, it is of interest for numerical analysis of non-linear chaotic behaviour. The equations are formulated as
\begin{align}
\frac{dX}{dt} &= -Y^2-Z^2-aX+aF, \label{eq:L1} \\
\frac{dY}{dt} &= XY-bXZ-Y+G, \nonumber \\
\frac{dZ}{dt} &= bXY+XZ-Z \nonumber
\end{align}
The variables $X(t)$, $Y(t)$, and $Z(t)$ represent certain meteorological systems such as wind currents and large-scale eddies. The coefficients $a$, $b$, $F$, and $G$ are chosen within certain limits in order to act as damping, coupling, and amplification of the physical processes. In Figure 2 an SGWRM solution for $a=0.25$, $b=4.0$, $F=8.0$ and $G=1.0$ is shown. Initial conditions are chosen as $(x,y,z)=(0.96,-1.1,0.5)$.

\begin{figure}[h!]
	\centering
	\includegraphics[width=5in]{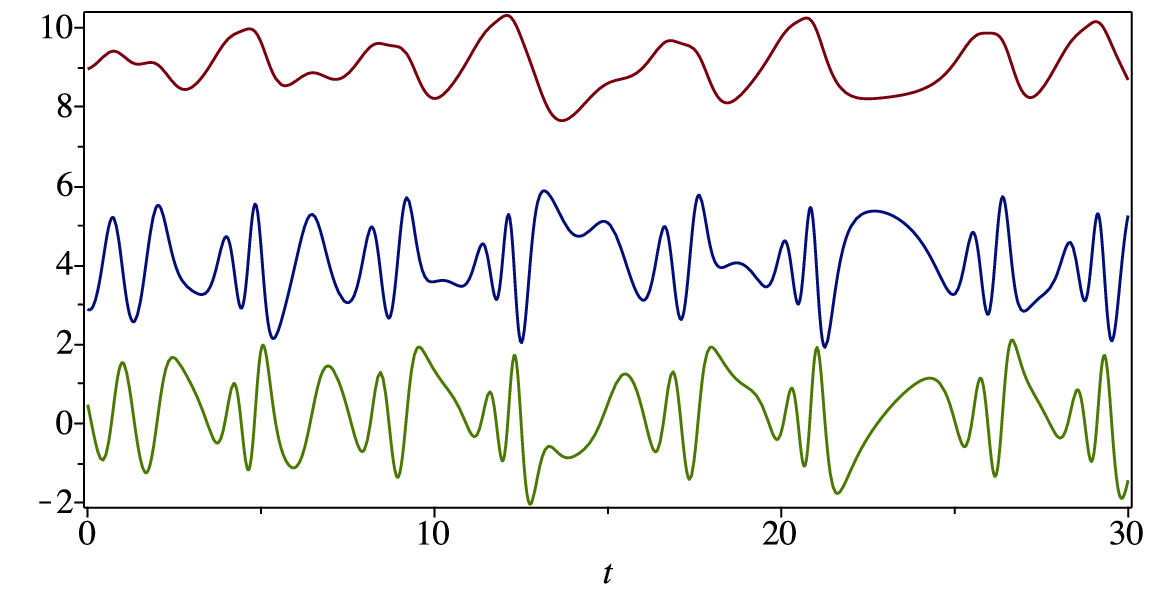}
	\caption{SGWRM solution of the Lorenz equations with parameters $a=0.25$, $b=4.0$, $F=8.0$ and $G=1.0$. The initial conditions are $(X,Y,Z)=(0.96,-1.1,0.5)$. From top to bottom; $X(t)+8$, $Y(t)+4$, $Z(t)$. The time-adaptive GWRM uses 61 time intervals and $K=8$ to obtain an accuracy $\epsilon=0.001$.}
\end{figure}

Explicit methods find these equations difficult to solve, since solution uncertainties at time $t$ will be strongly amplified by chaoticity, that is divergence from the true solution, in any interval $[t,t+\tau]$ unless the solution algorithm is very accurate. To see this, we make the following analysis, based on LLEs \cite{Scheffel_transforming}. We let $\delta u_1 \equiv u_1(t;\epsilon)-u_1(t;0) =  x(t;\epsilon)-x(t;0) \equiv \delta x$ and similarly for the $y$- and $z$- components. The following system of ODEs for $\bf \delta u$ is obtained;
\begin{equation}
	\begin{split}
	\frac{d \delta x}{dt} = -2y\delta y -2z\delta z - a\delta x  \\
	\frac{d \delta y}{dt} = x\delta y + y\delta x - bx\delta z - bz\delta x - \delta y  \\
	\frac{d \delta z}{dt} = bx\delta y + by\delta x + x \delta z + z\delta x - \delta z 
	\end{split}
\end{equation}

The coefficients for $\delta u_i$, including $x,y,z$, are assumed constant in each interval $[t,t+\tau]$ that is solved for. They do, of course, take on different values depending on the interval being investigated. 

The Lyapunov exponents can be determined as follows. Assume the linear system of ODE's

\begin{equation}
	\frac{d\bf u}{dt} = A\bf u
\end{equation}

\noindent where $\textbf u \in R^n$ and $\bf{A}$ is a constant, quadratic, diagonalizable matrix  with distinct eigenvalues ${\gamma}_k \in C$ and eigenvectors $\textbf{c}_k \in C^n$. It is easily shown that the solution to Eq. (25) is
\begin{equation}
	\textbf {u}(t) = \sum_{k=1}^na_k e^{\gamma_k} \textbf{c}_k
\end{equation}

\noindent in which $a_k$ are arbitrary constants. By comparison with Eq. (24) it is clear that $\textbf{u} = \delta \textbf{u}$ and $A=J$, where $J$ is the Jacobian of the system. The eigenvalues $\gamma_k$ are obtained from the characteristic equation 

\begin{equation}
	|J-\gamma I| = 0
\end{equation}

\noindent where $I$ is the identity matrix. Here

\begin{equation*}
J = 
\begin{pmatrix}
-a & -2y & -2z \\
y-bz & x-1 & -bx \\
by+z & bx & x-1
\end{pmatrix}
\end{equation*}

\noindent The Lorenz equations are chaotic at $t=0$ for these initial values; we obtain $\gamma_1=1.9$, $\gamma_2=-1.1+4.5i$ and $\gamma_3=-1.1-4.5i$. The parameters $F$ and $G$ do not affect the character of the ODE's, as seen from Eq. (24). Chaoticity is always manifested by at least one eigenvalue being positive for any value of $b$. Choosing $a \geq 3.1$, however, the system will initially be non-chaotic. In this case non-chaoticity also obtains for example for $t \approx 13-14$. The equations are not stiff for any values of $a$ and $b$.

In \cite{Scheffel:GWRM3} Eq. (23) is solved by time-adaptive GWRM, as well as by explicit RK4 and implicit Gauss-Legendre and Lobatto methods. The RK4 method is the most efficient of these. The GWRM however employs time intervals that are two orders of magnitude larger than the time steps of RK4, thus being up to four times faster at high accuracy. When computing perturbed scenarios, being of interest in numerical weather prediction, enhanced GWRM convergence further contributes to the superiority of the GWRM. A bonus feature is that the GWRM solutions have the form of analytically tractable Chebyshev series expansions.      

Chaoticity is consequentially problematic in numerical solution of ODEs and PDEs, not only because of the sensitivity to initial data, but also due to that local numerical errors become amplified with time. The question thus arises whether the implicitness of the GWRM, as mentioned in Section 2, mitigates numerical chaoticity or if transform methods are preferably applied \cite{Scheffel_transforming}. It was shown in \cite{Scheffel_transforming} that it is possible to, at least locally, \textit{transform the numerically generated chaos} of the equations to obtain aymptotically stable solutions in the transformed domain. This will be dealt with in the coming section.

\section{Achievable gain from time-spectral methods}

\noindent In an evaluation of the performance of time-spectral methods in general, and the GWRM in particular, all properties of solution steepness, solution oscillations as well as problem stiffness and chaoticity need to be taken into account. By "performance" we refer to accomplishment of convergence, accuracy and efficiency. We touch briefly upon some convergence and accuracy properties, whereafter we focus on efficiency for different types of problems.

\subsection{Convergence and accuracy}

\noindent Regarding \textit{convergence}, it was shown in \cite{Scheffel:GWRM7} that the SGWRM does not necessarily converge to an average over the sought solution in the solution interval. This is due to the fact that the time-spectral approximation uses a finite set of Chebyshev modes with collocation points at curve points, which generally do not represent the averaged solution. Furthermore, for strongly nonlinear problems the length of the time interval affects the convergence for iteratively determining the coefficients $a_{klm}$ in Eq. (9). Interval lengths may, however, be up to two order of magnitudes longer than those of explicit methods, even for chaotic problems \cite{Scheffel:GWRM3}. 

Concerning \textit{accuracy}, the GWRM employs Chebyshev polynomials which converge very rapidly: the error of a Chebyshev series for an infinitely differentiable function, truncated after $K$ terms, goes to zero more rapidly than any finite power of $1/K$ as $K\rightarrow \infty$ \cite{Gottlieb:1}. 

It is interesting to consider the built-in accuracy condition of the GWRM. Since the numerical values of Chebyshev polynomials are limited to the interval $[-1,1]$, the convergence requirement

\begin{equation}
	\frac{|a_{K-1,0,0}|+|a_{K,0,0}|}{|a_{0,0,0}|+|a_{1,0,0}|} < \epsilon ,
\end{equation}

\noindent routinely used for time-adaptive GWRM solution of PDEs, in turn corresponds to a solution relative accuracy. The solution error is typically up to an order of magnitude smaller, however, since, due to the strong Chebyshev series convergence, omitted terms $\ll \epsilon$. The accuracy of GWRM solutions are thus found without need to compute and compare with a solution produced with higher accuracy. This property is valuable for designing efficient time-adaptive algorithms; the prescribed accuracy can be directly implemented by adjusting the time interval length using the above criterion. Increased or reduced step lengths are conveniently set at little extra cost.

One may generally assume that for a prescribed GWRM accuracy, there is an essentially linear relationship between the number of required Chebyshev modes and the interval lengths. In other words, assuming that $K$ modes are required for $\epsilon$ accuracy in an interval of length $T$, a reduction of the interval length by some factor $k$ should approximately result in a corresponding reduction of the number of required modes. This is indeed found in numerical experiments, including the chaotic Lorenz equations solved in Section 4.2. To be more precise, only for intervals requiring $K \geq K^*$, where $K^* \approx 30$ for the Lorenz equations, the linear relationship is somewhat broken; for these long intervals additional modes are required.

\subsection{Efficiency}

\noindent In order to assess the \textit{efficiency} of the GWRM, as compared to finite difference methods, an effort has earlier been carried out to optimize the algorithms of the GWRM \cite{Scheffel:GWRM4}. Since the GWRM heavily rests on efficient iterative determination of the coefficients in Eq. (9), the nonlinear root solver SIR \cite{Scheffel:SIR,Scheffel:SIR2} was optimized. This entails careful selection of linear equation solvers, avoidance of costly re-computations of matrices after a few iterations, the use of band matrix methods, differentiation of non-zero band matrix elements only, use of spatial and temporal subdomains, and careful choice of initial iterate vector. As for the GWRM algorithm itself, primarily the use of several temporal and spatial subdomains with corresponding low Chebyshev order, overlapping rather than matched spatial subdomains, adaptive temporal subdomains, time parallelization and the use of Clenshaw's algorithm for Chebyshev series computations have proven fruitful. 

These are general measures for improving efficiency. It is now of interest to take a closer look at GWRM efficiency in relation to specific problem properties like solution steepness and oscillations as well as stiffness and chaoticity.

\subsection{GWRM solution of non-smooth problems}
\noindent Whereas time-spectral methods efficiently solves smooth problems, limited convergence for non-smooth problems with partially steep or highly oscillative solutions will add to iteration costs. It would thus be valuable to efficiently apply averaging techniques to the GWRM. For this reason, the LTA, TI and TA smoothing methods have been developed. The question now arises whether these represent improvements to the SGWRM.

\subsubsection{The Time-Integration (TI) Method}

\noindent The rationale for the TI method (see Section 3.2) is to introduce smoothing by solving for integrals of the sought solution, rather than solving for the solution itself, whereafter a back transformation through differentiation is carried out. It was found in \cite{Scheffel:GWRM7} that a moderate reduction in CPU time, as compared to SGWRM, could be obtained by a suitable choice of the smoothing parameter $A$ in the TI Method. This gain is due to the superior convergence properties of the TI Method, thus less time intervals are needed. As we will see, however, the performance of this method is not significantly superior to that of SGWRM. In Section 3 it was shown that the LTA and TI methods are essentially identical, since the computation of a long-time average (LTA method) also entails computing the time integral of the solution. The following general observations can be made. First, in order to obtain order $K$ accuracy in the solution, $u(t)$ in $du/dt$ must be expanded to order $K+1$. From the LTA Method integral Eq. (15), and the relation $u(t)=dZ/dt$ it is also seen $Z(t)$ must be expanded to order $K+1$ for the same accuracy. 

Similarly the function $v(t)$ of the TI Method, (Eq. (18)), must be expanded to order $K+1$ for this accuracy, since $u(t)=dv/dt-A$. It may be argued that the LTA and TI methods potentially can converge over larger time domains, but this is not necessarily helpful. Assume that for a particular case the TI Method converges in a single interval 2 times larger than those of the SGWRM, using an expansion order $K$ in each interval. Then approximately $2K$ modes are required by the TI Method for the same accuracy in the extended interval. For solving Eq. (9) using SIR, the CPU time scaling $N^{1.43}$, where $N$ is the number of modes, has been obtained \cite{Scheffel:GWRM4}. Thus the CPU time obtains as $(2K)^{1.43}=2.7K^{1.43}$ for the TI Method, in relation to $2K^{1.43}$ for the SGWRM. The latter relation uses the fact that two time intervals are solved for. 

CPU time gain is, however, not obtained indefinitely by reducing the interval length in order to reduce the number of required modes and iterations; the overhead cost per interval will take its due and an optimum interval length with respect to efficiency results. In summary; the TI Method cannot be superior to the SGWRM since there is no gain in CPU time from convergence at longer time intervals for the same accuracy. 

\subsubsection{The Time-Average (TA) Method}

\noindent Neither the TA Method is generally to be preferred to the SGWRM. This may seem surprising, since the sought, averaged, solution features a smoother dependence using less nodes and thus lower demands on the number of Chebyshev modes for a given accuracy. In \cite{Scheffel:GWRM7} some advantages of the TA Method were reported. These entail that for the chaotic Lorenz problem, the accuracy of the TA Method is less sensitive to initial conditions than SGWRM, and that high accuracy, as compared to the exact running average, is obtained even for strong averaging (large $\Delta$).

There is a problem, however, in that the number of equations to be solved in the TA Method are typically twice as many as those of SGWRM; thus also the number of Chebyshev modes of the algebraic system of coefficients Eq. (9) doubles. Furthermore the order of the combined system of differential equations may increase, demanding higher order Chebyshev expansions to retain accuracy. Neither is a combination of the TA Method and the TI Method favourable for increasing smoothness, since the the required number of Chebyshev modes for a certain accuracy within the TA Method cannot be reduced by the TI Method. In order to avoid doubling the number of differential equations, the LTA Method was suggested as an alternative averaging method to the TA method, but since the former is equivalent to the TI Method nothing is gained. 

As an efficient method for determining an approximate trend of the solution, however, without strong consideration of accuracy, the TA Method may be advantageous. It was seen, for example in \cite{Scheffel:GWRM7}, that only the TA Method could converge to the approximate (averaged) orbit of a charged particle moving in an inhomogeneous magnetic field for as long as three full (Larmor) orbits in a single time interval.

\subsection{Stiff and chaotic problems}

\noindent Time-spectral computation could have either of two purposes; 1) accurate determination of the temporal dependence in an interval, or 2) efficient approximate estimation of the trend of the solution. In the first case, we seek to minimize the required number of Chebyshev modes $K$ for a solution with prescribed relative error $\epsilon$ in an interval of length $T$. Minimization is essential, since the CPU time of the SIR equation solver, employed at each time interval, scales as $K^{1.43}$ and memory requirements scale as $K^{1.0}$ \cite{Scheffel:GWRM4}. The value of $K$ depends on the smoothness of the solution and on the number of extrema in the interval. In the second approach, accuracy $\epsilon$ may be sacrified in order to more efficiently follow the trend of the solution as it evolves in an initial-value problem. In this case we seek global SIR convergence for large $T$ at low $K$ values.

\subsubsection{GWRM solution of stiff problems}
\noindent Solving the stiff system of equations (21) it was seen that the SGWRM outperforms both standard explicit and implicit methods. This result is partly a manifestation of the high accuracy and long time intervals characteristic of time-spectral methods. But it is also a result of the implicit character of time-spectral methods, noted in Section 2. Time-spectral and implicit methods feature similar algoritmical structure, or stencils, as can be seen as follows. A generic 1D PDE in an unknown function $u$ has the form
\begin{equation}
	F(u,\frac{\partial u}{\partial t},\frac{\partial u}{\partial x},\frac{\partial^2 u}{\partial t \partial x},...) = 0 .
\end{equation}

Discretizing, let $u^n_i$ represent the value of the unknown function at grid point $i$ at time level $n$. Using a typical implicit scheme, the discretized equation may be formalized as
\begin{equation}
	G(u^{n+1}_i,u^{n+1}_{i-1},u^{n+1}_{i+1},u^{n}_i,u^{n}_{i-1},u^{n}_{i+1}) = 0  .
\end{equation}

This stencil includes not only spatially neighbouring points at the present time level $n$, but also points at the future time level $n+1$, at which the solution $u$ is solved for. The coefficients $a_{klm}$ of the GWRM solution ansatz Eq. (2) are obtained through iteration of Eq. (9). This entails finding the Chebyshev polynomial fit of the solution in the full time interval $[t_0,t_1]$, as can be seen from Eq. (4). Similarly as the function values $u^{n}_i,u^{n}_{i-1},u^{n}_{i+1}$ represents known initial values for the finite difference problem, the coefficients $b_{lm}$ in Eq. (9) represent initial values for the GWRM at $t_0$, from which $a_{klm}$, representing the solution in the interval, can be obtained iteratively. For the implicit finite difference scheme, a similar algorithm is used to obtain the updated $u$ values at point $n+1$.     

The implicit character of both methods symbolically requires inversion of system matrices and may limit algorithm efficiency. A number of optimization measures have thus been developed for the GWRM, as described in the beginning of this Section. Next, we consider the combined effect of GWRM implicitness, accuracy and efficiency for solving chaotic problems.

\subsubsection{GWRM solution of chaotic problems}
\noindent Chaotic systems in general, and the Lorenz equations of Section 4.2 in particular, are challenging with respect to numerical accuracy. It is well known that chaoticity typically requires small time steps in explicit methods due to the divergence of neighbouring solutions. 

In a recent study \cite{Morozova} it was furthermore found that, for implicit schemes, the presence of dense periodic phase space orbits in chaotic problems yields more stringent limits related to numerical accuracy than those related to numerical stability limits, so that all implicit schemes fail at approximately the same time-step sizes if these are a significant fraction of the quasi-period of the system. Unlike in laminar systems, only a modest gain in time step size is found for the implicit schemes when compared to the explicit fourth-order Runge-Kutta method. Computational cost per time-step largely determines the choice of an optimal implicit scheme. Among the implicit schemes being studied in \cite{Morozova}, one finds backward Euler, implicit Runge Kutta and trapezoidal methods.  

High numerical accuracy, or resolution, is a characteristic for time-spectral modelling of systems with smooth, or relatively smooth, solutions. The GWRM, for example, typically resolves a quasi-period of the Lorenz equation, at $\epsilon=0.001$ accuracy, using some 8 spectral modes, as shown in Section 3.3. Therefore, it becomes of interest to determine whether the GWRM is as robust in terms of stability as implicit schemes.  

From a computational point of view, the Lorenz system is third order in time, as can be seen by formally combining the equations. Employing the results of Section 3.3, 3 additional modes are required for Chebyshev expansion of each variable to obtain the same accuracy as for a direct Chebyshev approximation of the exact solution. From Fig. 2 it is seen that the $X(t),Y(t)$ and $Z(t)$ components of the exact solution feature 19, 44 and 41 extrema, respectively in the time interval [0,30]. Since 61 time intervals were employed, on average maximum $44/61=0.72$ extrema obtains per interval. Extrapolating Eq. (20), some $5.6+3=8.6$ modes are required for $\epsilon=0.001$ \textit{numerical resolution} accuracy of third order derivatives in each interval. In total, it is thus expected that $61 \times 8.6 \approx 520$ Chebyshev modes must be used for this accuracy. In the numerical computation of Figure 2, $K=8$ was used for each interval, so that $61 \times (8+1) \approx 550$ Chebyshev modes were employed. 

Considering that, as discussed in Section 5.1.3, the employed convergence criterion (28) is conservative in the sense that neglected terms $\ll \epsilon$, we may conclude that \textit{the SGWRM algorithm uses approximately as many modes that are required for numerical resolution only}; no additional modes are needed to resolve the chaotic behaviour of the differential equations.  

The reason for this behaviour is that the GWRM employs simultaneous (implicit) solution for the coefficients $a_{klm}$ of Eq. (2). This can be seen by rewriting Eq. (9) as
\begin{equation}
	a_{klm}=2\delta_{k0}b_{lm}+A_{klm}(\textbf{a})+F_{klm}
\end{equation}

\noindent where it holds that $(\textbf a)_{klm}=a_{klm}$. The coefficients $A_{klm}$ determine the spectral representation of the $D \textbf u$ term of Eq. (1). It is clearly seen that iteration of Eq. (31), using for example the SIR equation solver, corresponds to implicitly determining the solution for the entire temporal interval as defined by Eqs. (3) and (4). We can thus expect that, as shown here for the Lorenz equation, numerical resolution of the solution curve, rather than chaoticity, will determine the efficiency of the GWRM and other time-spectral methods. It was indeed found in \cite{Scheffel:GWRM3} that the SGWRM is more efficient than the implicit trapezoid and Lobatto IIIC methods, and the explicit RK4 method at $\epsilon < 0.001$, for solving the Lorenz problem (23).

\section{Discussion}

\noindent A major problem to handle when introducing exact smoothing methods for improving the accuracy and efficiency of time-spectral methods is the fact that the corresponding equations use as much information as the original equations. One can see this in two different ways. For the TA Method it must hold that 1) in order to calculate an exact average $U(t)$ according to the definition, there must be access to the exact form of $u(t)$ and 2) for small $\Delta$, $U(t)$ will approach $u(t)$, thus modelling all its details. This limits the efficiency of all the TI, LTA and TA methods. 

Common for these methods is that they are formally exact, which distinguish them from many previously suggested averaging algorithms.  A standard approach for weakly fluctuating solutions is to separate the solution into an averaged part and a fluctuating part which is eliminated ad hoc, analytically: $u(t)=U(t)+\tilde{u}(t)$, where $\tilde{u}(t)$ represents the fluctuating solution. By adding external information, the information to be computed becomes reduced, thus increasing efficiency. This method is, of course, of limited use for chaotic problems.

We can, in the following way, try to understand why smoothing methods, applied to the GWRM do not, as shown in Section 3, yield appreciable gain over SGWRM. Let us assume that an accuracy $\epsilon$ has been specified for the computed approximation $u^*(t)$ of the exact solution $u(t)$ of a differential equation in a sub-interval of length $T$. Assume that a direct Chebyshev approximation of the exact solution function requires $N$ modes for the given accuracy. Since the GWRM also entails a representation of the differential equation, a number of Chebyshev modes $N^*>N$ is required for $\epsilon$ accuracy; no GWRM method can be designed that would require $N$ modes or less. The potential gain in TI or TA methods is that a reduced number of modes $N^*_{opt}<N^*_{SGWRM}$ are sufficient, where $N^*_{SGWRM}$ denotes the required number of modes for SGWRM and $N^*_{opt}$ can be either of $N^*_{TI}$ or $N^*_{TA}$, or related to some other smoothing scheme. The corresponding gain in efficiency then amounts to $(N^*_{SGWRM}/N^*_{opt})^{\alpha}$, where $\alpha \in [1.4,3]$. Furthermore, convergence is of importance since if smaller time intervals need to be employed, then a large $N^*$ results. We can now make the following general statements. For $\epsilon \ll 1$ it is expected that $N^*_{opt}\rightarrow N^*_{SGWRM}$, using the TI Method. The TA Method should not be suitable in this limit. Accepting less accuracy, that is $\epsilon < 1$, it is not expected that $N^*_{opt}\rightarrow N$; rather the focus is on larger time intervals, for which the TI, TA or combinations of them, can satisfy the requirements of good convergence, low CPU cost and low memory consumption. The TI Method scaling (see Section 5.1.1) is however insufficient for obtaining higher efficiency than the SGWRM.

Furthermore, the TA Method potentially may reduce the number of extrema at high values of $\Delta$. In \cite{Scheffel:GWRM7} it was shown that the reduction of extrema, however, was only moderate and that this effect could not help to compensate for the increase in Chebyshev modes caused by the doubling of the number of equations to be solved in the TA Method. 

These results underscore the need for further exploration of alternative time-spectral approaches to enhance convergence and accuracy in the presence of steep or non-smooth solutions.

\section{Conclusion}

\noindent Time-spectral methods for solution of stiff and chaotic differential equations are compared to explicit and implicit finite difference methods with respect to efficiency. It is found that the time-spectral method GWRM demonstrates insensitivity to stiffness and chaoticity (diagnosed using local Lyapunov exponents) due to the implicit nature of its solution algorithm. Examples of efficient solution of these challenging problems, where explicit and implicit methods fail or are significantly slower, are given. Time-spectral methods are particularly advantageous when high accuracy is sought. 

Non-smooth and partially steep solutions, however, remain challenging for convergence and accuracy of time-spectral methods. Some, earlier suggested, smoothing algorithms are discussed and shown to be ineffective in addressing this issue. In the Time-Integration Method, the differential equations are reformulated as integrals of the sought solutions, which feature a smoother dependence than the sought solution itself. The Time-Average Method obtains smoothing by a reformulation to differential equations for the exact averages of the sought variables. Our findings underscore the need for further exploration of alternative time-spectral approaches to enhance convergence and accuracy in the presence of locally steep or non-smooth solutions.

\section{Data availability statement}
\noindent No datasets were generated or analysed during the current study.

\section{Declarations}
\noindent The author has no relevant financial or non-financial interests to disclose. \
\noindent No funding was received for conducting this study.

\section{References}

\label{}





\bibliographystyle{elsarticle-num}
\bibliography{mybib}







\end{document}